\documentclass{amsart}

\usepackage[frame,cmtip,arrow,matrix,line,graph,curve]{xy}
\usepackage{amsmath,amssymb, amscd}
\usepackage{amsfonts,amsxtra,amsthm}
\usepackage{mathrsfs}
\usepackage{eucal}
\usepackage{latexsym}

\numberwithin{equation}{section}

\newtheorem{theorem}{Theorem}[section]
\newtheorem{proposition}[theorem]{Proposition}
\newtheorem{corollary}[theorem]{Corollary}
\newtheorem{lemma}[theorem]{Lemma}

\newtheorem{remit}[theorem]{Remark}
\newtheorem{definit}[theorem]{Definition}

\newenvironment{remark}{\begin{remit}\rm}{\end{remit}}

\newcommand{\pp}{\mathbb{P}}
\newcommand{\qq}{\mathbb{Q}}
\newcommand{\cc}{\mathbb{C}}

\newcommand{\zz}{\mathbb{Z}}

\newcommand{\Hom}{\mathrm{Hom}}
\newcommand{\HHom}{\mathscr H\!om}

\newcommand{\Gr}{\mathrm{Gr}}
\newcommand{\Ext}{\mathrm{Ext}}
\newcommand{\EExt}{\mathscr E\!xt}

\newcommand{\git}{/\!\!/}

\newcommand{\cA}{\mathcal{A} }
\newcommand{\cB}{\mathcal{B} }
\newcommand{\cF}{\mathcal{F} }
\newcommand{\cE}{\mathcal{E} }

\newcommand{\cL}{\mathcal{L} }
\newcommand{\cO}{\mathcal{O} }

\newcommand{\tM}{\widetilde{M} }

\newcommand{\hM}{\widehat{M} }

\begin{document}

\title[Nonexistence of a crepant resolution of moduli spaces]
{Nonexistence of a crepant resolution of some moduli spaces of
sheaves on a K3 surface}
\date{}

\author{Jaeyoo Choy and Young-Hoon Kiem}
\address{Dept of Mathematics, Seoul National University,
Seoul 151-747, Korea} \email{donvosco@math.snu.ac.kr}
\email{kiem@math.snu.ac.kr}
\thanks{Young-Hoon Kiem was partially supported by KOSEF
R01-2003-000-11634-0; Jaeyoo Choy was partially supported by KRF
2003-070-C00001}

\keywords{Crepant resolution, irreducible symplectic variety,
moduli space, sheaf, K3 surface, desingularization, Hodge-Deligne
polynomial, Poincar\'{e} polynomial, stringy E-function}

\begin{abstract}
Let $M_c=M(2,0,c)$ be the moduli space of $\cO(1)$-semistable rank
2 torsion-free sheaves  with Chern classes $c_1=0$ and $c_2=c$ on
a K3 surface $X$ where $\cO(1)$ is a generic ample line bundle on
$X$. When $c=2n\geq4$ is even, $M_c$ is a singular projective
variety equipped with a holomorphic symplectic structure on the
smooth locus. In particular, $M_c$ has trivial canonical divisor.
 In \cite{ogrady}, O'Grady asks if
there is any symplectic desingularization of $M_{2n}$ for $n\ge
3$. In this paper, we show that there is no crepant resolution of
$M_{2n}$ for $n\geq 3$. This obviously implies that there is no
symplectic desingularization.
\end{abstract}
\maketitle

%%%%%%%%%%%%%%%%%%%%%%%%%%%%%%%%%%%%%%%%%%%%%%%%%%%%%%%%%%%%%%%%%%%
%%%%%%%%%%%%%%%%%%%%%%%%%%%%%%%%%%%%%%%%%%%%%%%%%%%%%%%%%%%%%%%%%%%

\section{Introduction}

Let $X$ be a complex projective K3 surface with polarization
$H=\cO_X(1)$ generic in the sense of \cite{ogrady} \S0. Let
$M(r,c_1,c_2)$ be the moduli space of rank $r$ $H$-semistable
torsion-free sheaves on $X$ with Chern classes $(c_1,c_2)$ in
$H^*(X,\zz)$. Let $M^s(r,c_1,c_2)$ be the open subscheme of
$H$-stable sheaves in $M(r,c_1,c_2)$. In \cite{Muk84}, Mukai shows
that $M^s(r,c_1,c_2)$ is smooth and has a holomorphic symplectic
structure. By \cite{Gi77}, if either $(c_1.H)$ or $c_2$ is an odd
number, then $M(2,c_1,c_2)$ is equal to $M^s(2,c_1,c_2) $ and thus
$M(2,c_1,c_2)$ is a smooth projective irreducible symplectic
variety. However if both $(c_1.H)$ and $c_2$ are even numbers then
generally $M(2,c_1,c_2)$ admits singularities. We restrict our
interest to the trivial determinant case i.e. $c_1=0$ and let
$M_c=M(2,0,c)$ where $c=2n$ ($n\geq2$). It is well-known that
$M_{2n}$ is an irreducible, normal (\cite{Yo01} Theorem 3.18) and
projective variety (\cite{HL97} Theorem 4.3.4) of dimension $8n-6$
(\cite{Muk84} Theorem 0.1) with only Gorenstein singularities
(\cite{HL97} Theorem 4.5.8, \cite{Ei95} Corollary 21.19). Since
$M_{2n}$ contains the smooth open subset $M^s_{2n}$, there arises
a natural question: does there exist a resolution of $M_{2n}$ such
that the Mukai form on $M^s_{2n}$ extends to the resolution
without degeneration? When $c=4$, O'Grady successfully extends the
Mukai form on $M^s_{2n}$ to some resolution without degeneration
(\cite{og97, ogrady}). At the same time, he conjectures
nonexistence of a symplectic desingularization of $M_{2n}$ for
$n\ge 3$ (\cite{ogrady}, (0.1)). Our main result in this paper is
the following.

\begin{theorem} \label{thm:main
result}  If $n\geq3$, there is no crepant resolution of $M_{2n}$.
\end{theorem}

The highest exterior power of a symplectic form gives a
non-vanishing section of the canonical sheaf on $M_{2n}$. Likewise
any symplectic desingularization of $M_{2n}$ has trivial canonical
divisor and hence it must be a crepant resolution. Therefore,
O'Grady's conjecture is a consequence of Theorem \ref{thm:main
result}.

\begin{corollary} \label{cor:O'Grady's conjecture} If $n\geq3$,
there is no symplectic desingularization of $M_{2n}$.
\end{corollary}

The idea of the proof of Theorem \ref{thm:main result} is to use a
new invariant called the stringy E-function \cite{Bat98, DL99}.
Since $M_{2n}$ is normal irreducible variety  with log terminal
singularities (\cite{ogrady}, 6.1), the stringy E-function of
$M_{2n}$ is a well-defined rational function. If there is a
crepant resolution $\tM_{2n}$ of $M_{2n}$, then the stringy
E-function of $M_{2n}$ is equal to the Hodge-Deligne polynomial
(E-polynomial) of $\tM_{2n}$ (Theorem \ref{thm:Batyrev's result}).
In particular, we deduce that the stringy E-function
$E_{st}(M_{2n};u,v)$ must be a polynomial. Therefore, Theorem
\ref{thm:main result} is a consequence of the following.
\begin{proposition}\label{prop:stringy E-function test} The stringy
E-function $E_{st}(M_{2n};u,v)$ is not a polynomial for $n\geq3$.
\end{proposition} To prove that
$E_{st}(M_{2n};u,v)$ is not a polynomial for $n\geq 3$, we show
that $E_{st}(M_{2n};z,z)$ is not a polynomial in $z$. Thanks to
the detailed analysis of Kirwan's desingularization in \cite{og97}
and \cite{ogrady} which is reviewed in section 4, we can find an
expression for $E_{st}(M_{2n};z,z)$ and then with some efforts on
the combinatorics of rational functions we show that
$E_{st}(M_{2n};z,z)$ is not a polynomial in section 3. In section
2, we recall basic facts on stringy E-function and in section 5 we
prove a lemma which computes the E-polynomial of a divisor.

In \cite{ogrady}, O'Grady  gets a symplectic desingularization
$\tM_{2n}$ of $M_{2n}$ in the case when $n=2$. This turns out to
be a new irreducible symplectic variety, which means that it does
not come from a generalized Kummer variety nor from a Hilbert
scheme parameterizing 0-dimensional subschemes on a K3 surface
\cite{og98, Bea83}. Corollary \ref{cor:O'Grady's conjecture} shows
that unfortunately we cannot find any more irreducible symplectic
variety in this way.

After we finished the first draft of this paper, we learned that
Kaledin and Lehn \cite{KL04} proved Corollary \ref{cor:O'Grady's
conjecture} in a completely different way. We are grateful to D.
Kaledin for informing us of their approach. The second named
author thanks Professor Jun Li for useful discussions concerning
the article \cite{VW94}. Finally we would like to express our
gratitude to the referee for careful reading and challenging us
for many details which led us to improve the manuscript and
correct an error in Proposition 3.2.

\section{Preliminaries}
In this section we collect some facts that we shall use later.

For a topological space $V$, the Poincar\'e polynomial of $V$ is
defined as
\begin{equation} \label{eqn:Poincare polynomial}
P(V;z)=\sum_{i}(-1)^ib_i(V)z^i
\end{equation} where $b_i(V)$ is the $i$-th Betti number of $V$.
 It is well-known from \cite{Go90} that the Betti numbers of the
Hilbert scheme of points $X^{[n]}$ in $X$ are given by the
following:
\begin{equation}\label{eqn:Betti for X[n]} \sum_{n\geq
0}P(X^{[n]};z)t^n=\prod_{k\geq1}
\prod_{i=0}^{4}(1-z^{2k-2+i}t^k)^{(-1)^{i+1}b_i(X)}.
\end{equation}
% In particular, the Euler numbers $a_n$ of $X^{[n]}$
%are given by the equation
%\[\sum^\infty_{n=0}e_nt^n=\prod^\infty_{k=1}1/(1-t^k)^{24}.\]

Next we recall the definition and basic facts about stringy
E-functions from \cite{Bat98,DL99}. Let $W$ be a normal
irreducible variety with at worst log-terminal singularities, i.e.
\begin{enumerate} \item W is $\qq$-Gorenstein;
\item for a resolution of singularities $\rho: V\to W$ such that
the exceptional locus of $\rho$ is a divisor $D$ whose irreducible
components $D_1,\cdots,D_r$ are smooth divisors with only normal
crossings, we have \[K_V=\rho^*K_W+\sum^r_{i=1} a_iD_i \] with
$a_i>-1$ for all $i$, where $D_i$ runs over all irreducible
components of $D$. The divisor $\sum^r_{i=1}a_iD_i$ is called the
\textit{discrepancy divisor}. \end{enumerate}

For each subset $J\subset I=\{1,2,\cdots,r\}$, define
$D_J=\cap_{j\in J}D_j$, $D_\emptyset=V$ and $D^0_J=D_J-\cup_{i\in
I-J}D_i$. Then the stringy E-function of $W$ is defined by
\begin{equation} \label{eqn:stringy E-function}
E_{st}(W;u,v)=\sum_{J\subset I}E(D^0_J;u,v)\prod_{j\in
J}\frac{uv-1}{(uv)^{a_j+1}-1} \end{equation} where \[ E(Z;u,v) =
\sum_{p,q}\sum_{k\geq 0} (-1)^kh^{p,q}(H^k_c(Z;\cc))u^pv^q \] is
the Hodge-Deligne polynomial for a variety $Z$. Note that the
Hodge-Deligne polynomials have \begin{enumerate} \item the
additive property: $E(Z;u,v)=E(U;u,v)+E(Z-U;u,v)$ if $U$ is a
smooth open subvariety of $Z$; \item the multiplicative property:
$E(Z;u,v)=E(B;u,v)E(F;u,v)$ if $Z$ is a Zariski locally trivial
$F$-bundle over $B$. \end{enumerate}

By \cite{Bat98} Theorem 6.27, the function $E_{st}$ is independent
of the choice of a resolution (Theorem 3.4 in \cite{Bat98}) and
the following holds.
\begin{theorem} \label{thm:Batyrev's result} (\cite{Bat98} Theorem
3.12) Suppose $W$ is a $\qq$-Gorenstein algebraic variety with at
worst log-terminal singularities. If $\rho:V\to W$ is a crepant
desingularization (i.e. $\rho^*K_W=K_V$) then
$E_{st}(W;u,v)=E(V;u,v)$. In particular, $E_{st}(W;u,v)$ is a
polynomial.
\end{theorem}

\section{Proof of Proposition \ref{prop:stringy E-function test}}
In this section we first find an expression for the stringy
E-function of the moduli space $M_{2n}$ for $n\geq 3$ by using the
detailed analysis of Kirwan's desingularization in \cite{og97,
ogrady}. Then we show that it cannot be a polynomial, which proves
Proposition \ref{prop:stringy E-function test}.

We fix a generic polarization of $X$ as in \cite{ogrady}. The
moduli space $M_{2n}$ has a stratification \[M_{2n}=M^s_{2n}\sqcup
(\Sigma-\Omega)\sqcup \Omega\] where $M^s_{2n}$ is the locus of
stable sheaves and $\Sigma\simeq(X^{[n]}\times X^{[n]})/{\rm
involution}$ is the locus of sheaves of the form $I_Z\oplus
I_{Z'}$ ($[Z],[Z']\in X^{[n]}$) while $\Omega\simeq X^{[n]}$ is
the locus of sheaves $I_Z\oplus I_{Z}$. For $n\geq 3$, Kirwan's
desingularization $\rho:\hM_{2n}\to M_{2n}$ is obtained by blowing
up $M_{2n}$ first along $\Omega$, next along the proper transform
of $\Sigma$ and finally along the proper transform of a subvariety
$\Delta$ in the exceptional divisor of the first blow-up. This is
indeed a desingularization by \cite{ogrady} Proposition 1.8.3.

Let $D_1=\hat{\Omega}$, $D_2=\hat{\Sigma}$ and $D_3=\hat{\Delta}$
be the (proper transforms of the) exceptional divisors of the
three blow-ups. Then they are smooth divisors with only normal
crossings as we will see in Proposition \ref{prop:analysis on exc}
and the discrepancy divisor of $\rho:\hM_{2n}\to M_{2n}$ is
(\cite{ogrady}, 6.1)
\[(6n-7)D_1+(2n-4)D_2+(4n-6)D_3.
\] Therefore the singularities are log-terminal for $n\geq 2$, and
from (\ref{eqn:stringy E-function}) the stringy E-function of
$M_{2n}$ is given by \begin{eqnarray}\label{eqn:stringy E-function
of M_c}
E(M^s_{2n};u,v)+E(D^0_1;u,v){\textstyle\frac{1-uv}{1-(uv)^{6n-6}}}
+E(D^0_2;u,v){\textstyle\frac{1-uv}{1-(uv)^{2n-3}}}\nonumber \\
+E(D^0_3;u,v) {\textstyle\frac{1-uv}{1-(uv)^{4n-5}}}
+E(D^0_{12};u,v){\textstyle\frac{1-uv}{1-(uv)^{6n-6}}\frac{1-uv}{1-(uv)^{2n-3}}} \\
+E(D^0_{23};u,v){\textstyle\frac{1-uv}{1-(uv)^{2n-3}}\frac{1-uv}{1-(uv)^{4n-5}}}
+E(D^0_{13};u,v){\textstyle\frac{1-uv}{1-(uv)^{4n-5}}\frac{1-uv}{1-(uv)^{6n-6}}}
\nonumber \\
+E(D^0_{123};u,v){\textstyle\frac{1-uv}{1-(uv)^{6n-6}}
\frac{1-uv}{1-(uv)^{2n-3}}\frac{1-uv}{1-(uv)^{4n-5}}} .\nonumber
\end{eqnarray}
% and from (\ref{eqn:stringy Euler number}) the stringy
%Euler number of $M_{2n}$ is given by
%\begin{eqnarray}\label{eqn:stringy Euler number for M_c}
%e(M^s_{2n})+e(D^0_1){\textstyle\frac1{6n-6}}
%+e(D^0_2){\textstyle\frac1{2n-3}}+e(D^0_3){\textstyle\frac1{4n-5}} \nonumber \\
%+e(D^0_{12}){\textstyle\frac1{6n-6}\frac1{2n-3}}
%+e(D^0_{23}){\textstyle\frac1{2n-3}\frac1{4n-5}} \\
%+e(D^0_{13}){\textstyle\frac1{4n-5}\frac1{6n-6}}
%+e(D^0_{123}){\textstyle\frac1{6n-6}\frac1{2n-3}\frac1{4n-5}}.\nonumber
%\end{eqnarray}

We need to compute the Hodge-Deligne polynomials of $D^0_J$ for
$J\subset \{1,2,3\}$. Let $(\cc^{2n},\omega)$ be a symplectic
vector space. Let $\Gr^{\omega}(k,2n)$ be the Grassmannian of
$k$-dimensional subspaces of $\cc^{2n}$, isotropic with respect to
the symplectic form $\omega$ (i.e. the restriction of $\omega$ to
the subspace is zero).

\begin{lemma}\label{lem:Hodge poly of Gr} For $k\leq n$, the Hodge-Deligne polynomial of
$\Gr^\omega(k,2n)$ is
\[\prod_{1\leq i\leq k} \frac{1-(uv)^{2n-2k+2i}}{1-(uv)^i}. \]
%Hence the Euler number of $\Gr^\omega(k,2n)$ is $2^k {n\choose
%k}$.
\end{lemma}

\proof Consider the incidence variety \[ Z= \{(a,b)\in
\Gr^\omega(k-1,2n)\times \Gr^\omega(k,2n)|a\subset b\}. \] This is
a $\pp^{2n-2k+1}$-bundle over $\Gr^\omega(k-1,2n)$ and a
$\pp^{k-1}$-bundle over $\Gr^\omega(k,2n)$. We have the following
equalities between Hodge-Deligne polynomials: \begin{eqnarray*}
E(Z;u,v)&=&\frac{1-(uv)^{2n-2k+2}}{1-uv}
E(\Gr^\omega(k-1,2n);u,v)\\ &=& \frac{1-(uv)^{k}}{1-uv}
E(\Gr^\omega(k,2n);u,v). \end{eqnarray*} The desired formula
follows recursively from $\Gr^\omega(1,2n)=\pp^{2n-1}$. \qed\\
%
%In particular, the Poincar\'{e} polynomial of  $\Gr^\omega(k,2n)$
%is
%\begin{equation} \label{eqn:equality of Poincare polynomial}
%P(\Gr^\omega(k,2n);z)= \prod_{1\leq i\leq k}
%\frac{1-z^{2(2n-2k+2i)}}{1-z^{2i}}.
%\end{equation}

Let $\hat{\pp}^5$ be the blow-up of $\pp^5$ (projectivization of
the space of $3\times 3$ symmetric matrices) along $\pp^2$ (the
locus of rank 1 matrices). We have the following from \cite{og97}
and \cite{ogrady}. The proof will be presented in \S\ref{sec:
Proof of Lemma}.

\begin{proposition}\label{prop:analysis on exc}
Let $n\geq 3$.

(1) $D_1$ is a $\hat{\pp}^5$-bundle over a
$\Gr^\omega(3,2n)$-bundle over $X^{[n]}$.

(2) $D_2^0$ is a free $\zz_2$-quotient of a Zariski locally
trivial $I_{2n-3}$-bundle over $ X^{[n]}\times
X^{[n]}-\mathbf{\Delta} $ where $\mathbf{\Delta}$ is the diagonal
in $ X^{[n]}\times X^{[n]}$ and $I_{2n-3}$ is the incidence
variety given by
\[ I_{2n-3}=\{(p,H)\in \pp^{2n-3}\times \breve{\pp}^{2n-3}| p\in
H\}.
\]

(3) $D_3$ is a $\pp^{2n-4}$-bundle over a Zariski locally trivial
$ \pp^2$-bundle over a Zariski locally trivial
$\Gr^\omega(2,2n)$-bundle over $X^{[n]}$.

(4) $D_{12}$ is a $\pp^2$-bundle over a $\pp^2$-bundle over a
$\Gr^\omega(3,2n)$-bundle over $X^{[n]}$.

(5) $D_{23}$ is a $\pp^{2n-4}$-bundle over a $ \pp^1$-bundle over
a $\Gr^\omega(2,2n)$-bundle over $X^{[n]}$.

(6) $D_{13}$ is a $ \pp^2$-bundle over a $\pp^2$-bundle over a
$\Gr^\omega(3,2n)$-bundle over $X^{[n]}$.

(7) $D_{123}$ is a $\pp^1$-bundle over a $\pp^2$-bundle over
a $\Gr^\omega(3,2n)$-bundle over $X^{[n]}$. \\
All the above bundles except in (2) and (3) are Zariski locally
trivial. Moreover, $D_i$ ($i=1,2,3$) are smooth divisors such that
$D_1\cup D_2\cup D_3$ is normal crossing.
\end{proposition}

From Lemma \ref{lem:Hodge poly of Gr} and Proposition
\ref{prop:analysis on exc}, we have the following corollary by the
additive and multiplicative properties of the Hodge-Deligne
polynomial.

\begin{corollary}\label{eqn:computation of stringy E-function}
$$ E(D_1;u,v) = \Bigl({\textstyle
\frac{1-(uv)^6}{1-uv}-\!\frac{1-(uv)^3}{1-uv}+\!\bigl(\frac{1-(uv)^3}{1-uv}\bigr)^2}\Bigr)
\! \times\!\!\! \prod_{1\leq i\leq 3}\! \Bigl({\textstyle
\frac{1-(uv)^{2n-6+2i}}{1-(uv)^i}}\Bigr)\!\! \times\!
E(X^{[n]};u,v),$$

$$E(D_3;u,v)   = {\textstyle
\frac{1-(uv)^{2n-3}}{1-uv}\cdot\frac{1-(uv)^3}{1-uv}} \times
\prod_{1\leq i\leq 2}\Bigl({\textstyle \frac{1-(uv)^{2n-4+2i}}
{1-(uv)^i}}\Bigr)\times E(X^{[n]};u,v), $$

$$ E(D_{12};u,v)   = \Bigl({\textstyle
\frac{1-(uv)^3}{1-uv}}\Bigr)^2\times \prod_{1\leq i\leq
3}\Bigl({\textstyle \frac{1-(uv)^{2n-6+2i}}{1-(uv)^i}}\Bigr)
\times E(X^{[n]};u,v), $$

$$ E(D_{23};u,v)  = {\textstyle
\frac{1-(uv)^{2n-3}}{1-uv}\cdot\frac{1-(uv)^2}{1-uv}}
\times\prod_{1\leq i\leq 2}\Bigl({\textstyle
\frac{1-(uv)^{2n-4+2i}}{1-(uv)^i}}\Bigr)\times E(X^{[n]};u,v), $$

$$E(D_{13};u,v)   = {\textstyle \frac{
1-(uv)^3}{1-uv}\cdot\frac{1-(uv)^{2n-4}}{1-uv}} \times
\prod_{1\leq i\leq 2}\Bigl({\textstyle
\frac{1-(uv)^{2n-4+2i}}{1-(uv)^i}}\Bigr)\times E(X^{[n]};u,v), $$
$$ E(D^0_{123};u,v)   ={\textstyle
\frac{1-(uv)^2}{1-uv}\cdot\frac{1-(uv)^{2n-4}}{1-uv}}\times
\prod_{1\leq i\leq 2}\Bigl({\textstyle
\frac{1-(uv)^{2n-4+2i}}{1-(uv)^i}}\Bigr)\times E(X^{[n]};u,v).$$
\end{corollary}
\proof Perhaps the only part that requires proof is the equation
for $E(D_3;u,v)$. From Proposition \ref{prop:analysis on exc} (3),
$D_3$ is a projective variety which is a $\pp^{2n-4}$-bundle over
a smooth projective variety, say $Y$, whose E-polynomial is
$$E(\pp^2;u,v)\times E(\mathrm{Gr}^\omega (2,2n);u,v)\times
E(X^{[n]};u,v).$$ By the Leray-Hirsch theorem (\cite{V02I} p.182),
we have
\begin{eqnarray*}H^*(D_3;\cc)\cong H^*(Y;\cc)\otimes
H^*(\pp^{2n-4};\cc)\cong H^*(Y;\cc)\otimes \cc
[\lambda]/(\lambda^{2n-3})\\ \cong H^*(Y;\cc)\oplus
H^*(Y;\cc)\lambda\oplus \cdots \oplus
H^*(Y;\cc)\lambda^{2n-4}\end{eqnarray*} where $\lambda$ is a class
of type $(1,1)$ which comes from the K\"ahler class. The above
determines the Hodge structure of $D_3$ because the Hodge
structure is compatible with the cup product. Therefore we deduce
that $$ E(D_3;u,v)   = {\textstyle
\frac{1-(uv)^{2n-3}}{1-uv}\times E(Y;u,v)}.  $$ \qed

For the E-polynomial of $D_2^0$ we have the following lemma whose
proof is presented in section \ref{sec: Computation of E-poly of
D_0^2}. Recall that $$I_{2n-3}=\{((x_i),(y_j))\in \pp^{2n-3}\times
\pp^{2n-3}\,|\, \sum_{i=0}^{2n-3} x_iy_i=0\}$$ and there is an
action of $\zz_2$ which interchanges $(x_i)$ and $(y_j)$. Let
$H^r(I_{2n-3})^+$ denote  the $\zz_2$-invariant subspace of
$H^r(I_{2n-3})$ .
\begin{lemma}\label{lem: Hodge Deligne poly of D02}
\begin{eqnarray}\label{eqn: compute D02} \lefteqn{
 E(D^0_2;z,z)=P(I_{2n-3};z) \Bigl(
 \frac{P(X^{[n]};z)^2-P(X^{[n]};z^2)}2 \Bigr)} && \\ && +
 P^+(I_{2n-3};z)\bigl(P(X^{[n]};z^2)-P(X^{[n]};z) \bigr)\nonumber
 \end{eqnarray}
where  $P^+(I_{2n-3};z)=\displaystyle \sum_{r\geq0}(-1)^rz^r\dim
H^r(I_{2n-3})^+$.  Moreover
\begin{eqnarray}\label{eqn: E D02 is divisible by some Q}
 E(D^0_2;z,z)=\frac{1-(z^2)^{2n-3}}{1-z^2} Q(z^2)\end{eqnarray} for some
polynomial $Q$.
\end{lemma}

\textit{Proof of Proposition \ref{prop:stringy E-function test}.}

Let us prove that \eqref{eqn:stringy E-function of M_c} cannot be
a polynomial. Let
$$S(z)=E_{st}(M_{2n};z,z)-E(M^s_{2n};z,z).$$ It suffices to show
that $S(z)$ is not a polynomial for all $n\geq3$ because
$E(M^s_{2n};z,z)$ is a polynomial.

Note that given any $n\geq 3$, we can explicitly compute
$E(X^{[n]};z,z)$ and $E(D^0_2;z,z)$ by (\ref{eqn:Betti for X[n]})
and Lemma \ref{lem: Hodge Deligne poly of D02}. If $n=3$, direct
calculation shows that $S(z)$ is as follows:
\begin{eqnarray*} S(z)& =&
 1+46z^2+852z^4+12308z^6+111641z^8+886629z^{10}+4233151z^{12}\\
 & & +4990239z^{14}+4999261z^{16}+4230852z^{18}+884441z^{20}+113877z^{22}\\
 & & +12928z^{24}+3749z^{26}+3200z^{28}+2877z^{30}+299z^{32}+\cdots.
\end{eqnarray*} It is easy to see from (\ref{eqn:stringy E-function of M_c})
and Corollary \ref{eqn:computation of stringy E-function} that if
$S(z)$ were a polynomial, it should be of degree $\le 30$. Since
the series $S(z)$ has a nonzero coefficient of $z^{32}$, $S(z)$
cannot be a polynomial. So we assume from now on that $n\ge 4$.

Express the rational function $S(z)$ as
$$\frac{N(z)}{(1-(z^2)^{2n-3})(1-(z^2)^{4n-5})(1-(z^2)^{6n-6})}.$$
All we need to show is that the numerator $N(z)$ is not divisible
by the denominator
$(1-(z^2)^{2n-3})(1-(z^2)^{4n-5})(1-(z^2)^{6n-6})$.

As $E(X^{[n]};z,z)$ and $E(D^0_2;z,z)$ do not have nonzero terms
of odd degree by (\ref{eqn:Betti for X[n]}) and Lemma \ref{lem:
Hodge Deligne poly of D02}, all the nonzero terms in $S(z)$ are of
even degree by (\ref{eqn:stringy E-function of M_c}) and Corollary
\ref{eqn:computation of stringy E-function}. Hence, we can write
$S(z)=s(z^2)=s(t)$ for some rational function $s(t)$ in $t=z^2$.
The numerator $N(z)=n(z^2)=n(t)$ is not divisible by
$1-(z^2)^{2n-3}$ if and only if $n(t)$ is not divisible by
$1-t^{2n-3}$. By direct computation using (\ref{eqn:stringy
E-function of M_c}), Corollary \ref{eqn:computation of stringy
E-function} and Lemma \ref{lem: Hodge Deligne poly of D02}, $n(t)$
modulo $1-t^{2n-3}$ is congruent to
\begin{eqnarray}\label{eqn:denumerator modulo}
\shoveleft(1-t)^2(1-t^{4n-5})\times\Bigl({\textstyle
\frac{1-t^3}{1-t}}\Bigr)^2\times \prod_{1\leq i\leq
3}\Bigl({\textstyle \frac{1-t^{2n-6+2i}}{1-t^i}}\Bigr) \times
p(X^{[n]};t) \\ -(1-t)^2(1-t^{4n-5})\times {\textstyle
\frac{1-t^2}{1-t}\cdot\frac{1-t^{2n-4}}{1-t}}\times \prod_{1\leq
i\leq 2}\Bigl({\textstyle \frac{1-t^{2n-4+2i}}{1-t^i}}\Bigr)
\times p(X^{[n]};t) \nonumber  \\ -(1-t)^2(1-t^{6n-6})\times
{\textstyle \frac{1-t^2}{1-t}\cdot\frac{1-t^{2n-4}}{1-t}}\times
\prod_{1\leq i\leq 2}\Bigl({\textstyle
\frac{1-t^{2n-4+2i}}{1-t^i}}\Bigr) \times p(X^{[n]};t) \nonumber
\\ + (1-t)^3\times{\textstyle \frac{1-t^2}{1-t}
\cdot\frac{1-t^{2n-4}}{1-t}}\times \prod_{1\leq i\leq
2}\Bigl({\textstyle \frac{1-t^{2n-4+2i}}{1-t^i}}\Bigr)\times
p(X^{[n]};t) \nonumber  \end{eqnarray} where
$p(X^{[n]};t)=P(X^{[n]};z)$ with $t=z^2$. We write
(\ref{eqn:denumerator modulo}) as a product $\bar s(t)\cdot
p(X^{[n]};t)$ for some polynomial $\bar s(t)$. For the proof of
our claim for $n\geq 4$, it suffices to prove the following:
\begin{enumerate} \item if $n$ is not divisible by 3, then $1-t$
is the GCD of $1-t^{2n-3}$ and $\bar s(t)$, and
$\frac{1-t^{2n-3}}{1-t}$ does not divide $p(X^{[n]};t)$; \item if
$n$ is divisible by 3, then $1-t^3$ is the GCD of $1-t^{2n-3}$ and
$\bar s(t)$, and $\frac{1-t^{2n-3}}{1-t^3}$ does not divide
$p(X^{[n]};t)$.
\end{enumerate}

For (1), suppose $n$ is not divisible by 3. From
(\ref{eqn:denumerator modulo}), $\bar s(t)$ is divisible by $1-t$.
We claim that $\bar s(t)$ is not divisible by any irreducible
factor of $\frac{1-t^{2n-3}}{1-t}$, i.e. for any root $\alpha$ of
$1-t^{2n-3}$ which is not 1, $\bar s(\alpha)\neq 0$. Using the
relation $\alpha^{2n-3}=1$, we compute directly that
\begin{equation}\label{eqn:bar s} \bar s(\alpha)={\textstyle
-\frac{\alpha(1-\alpha^{-1}){(1-\alpha^3)}^2}{1+\alpha}},
\end{equation} which is not 0 because 3 does not divide $2n-3$.

Next we check that $\frac{1-t^{2n-3}}{1-t}$ does not divide
$p(X^{[n]};t)$. %Note that $p(X^{[n]};t)$ is a palindromic
%polynomial of degree $2n$ by Poincar\'{e} duality.
We put
$${\displaystyle p(X^{[n]};t)=\sum_{0\leq i\leq 2n} c_it^i}$$ and
write $p(X^{[n]};t)$ as follows:
\begin{eqnarray} \label{eqn:p(t) when 3 NOT divides n} \lefteqn{
\sum_{0\leq i\leq 2n} c_it^i  =
(c_0+c_{2n-3})+(c_1+c_{2n-2})t+(c_2+c_{2n-1})t^2 +(c_3+c_{2n})t^3
} \\ & & +  \sum_{4\leq i\leq 2n-4} c_it^i + c_{2n-3}(t^{2n-3}-1)
+  c_{2n-2}t(t^{2n-3}-1)\nonumber \\ & & +
c_{2n-1}t^2(t^{2n-3}-1)+ c_{2n}t^3(t^{2n-3}-1).\nonumber
\end{eqnarray} Therefore, the divisibility
of $p(X^{[n]};t)$ by $\frac{1-t^{2n-3}}{1-t}$ is that of
$(c_0+c_{2n-3})+ (c_1+c_{2n-2})t+ (c_2+c_{2n-1})t^2
+(c_3+c_{2n})t^3 + {\displaystyle \sum_{4\leq i\leq 2n-4} c_it^i}$
by $\frac{1-t^{2n-3}}{1-t}$. Since
$\frac{1-t^{2n-3}}{1-t}={\displaystyle \sum_{0\leq i\leq 2n-4}
t^i}$, the polynomial $(c_0+c_{2n-3})+ (c_1+c_{2n-2})t+
(c_2+c_{2n-1})t^2 +(c_3+c_{2n})t^3 + {\displaystyle \sum_{4\leq
i\leq 2n-4} c_it^i}$ is divisible by $\frac{1-t^{2n-3}}{1-t}$ if
and only if it is a scalar multiple of ${\displaystyle \sum_{0\leq
i\leq 2n-4} t^i}$, i.e. $c_0+c_{2n-3}=c_1+c_{2n-2}=
c_2+c_{2n-1}=c_3+c_{2n}= c_4=\cdots =c_{2n-4}$ ($n\geq 4$).

Table \ref{table:list of ci} is the list of $c_i$ ($1\leq i\leq
4$) for $n\geq 3$, which comes from direct computation using the
generating functions (\ref{eqn:Betti for X[n]}) for the Betti
numbers of $X^{[n]}$. By Table \ref{table:list of ci}, we can
check that this is impossible.  Indeed, for $n\geq 6$, $c_0=1$,
$c_1=23$, $c_2=300$ and $c_3=2876$, which implies $c_{2n-3}=2876$,
$c_{2n-2}=300$, $c_{2n-1}=23$ and $c_{2n-2}=1$ by Poincar\'{e}
duality. Thus $c_0+c_{2n-3}=2877$ while $c_1+c_{2n-2}=323$. For
$4\leq n\leq 5$, the proof is also direct computation using Table
\ref{table:list of ci}.

\begin{table}[t]
\begin{tabular}{c|cccccc}
      &$n=3$& $n=4$ & $n=5$& $n=6$& $n=7$&$n\geq 8$ \\ \hline
$c_1$ & 23  & 23    & 23   & 23   & 23   & 23 \\
$c_2$ & 299 & 300   & 300  &300   &300   &300 \\
$c_3$ & 2554& 2852  & 2875 & 2876 &2876  &2876 \\
$c_4$ & 299 & 19298 &22127 &22426 &22449 &22450
\end{tabular}\caption{list of $c_i$ \label{table:list of ci}}
\end{table}

For (2), suppose 3 divides $n$ and $n\neq 3$. Then from
(\ref{eqn:bar s}), $(1-t^3)$ divides $\bar s(t)$. More precisely,
for a third root of unity $\alpha$, $\bar s(\alpha)=0$. On the
other hand, if $\alpha$ is a root of $1-t^{2n-3}$ but not a third
root of unity then we can observe that $\bar s(\alpha)\neq 0$ by
(\ref{eqn:bar s}). Therefore, since every root of $1-t^{2n-3}$ is
a simple root, any irreducible factor of
$\frac{1-t^{2n-3}}{1-t^3}$ does not divide $\bar s(t)$.

We next check that the polynomial $\frac{1-t^{2n-3}}{1-t^3}$ does
not divide $p(X^{[n]};t)$. Write $p(X^{[n]};t)=\displaystyle
\sum_{0\leq i\leq 2n} c_it^i$ as follows:
\begin{eqnarray}\label{eqn:p(t) when 3 divides n} \lefteqn{\sum_{0\leq i\leq 2n} c_it^i   =
(c_0+c_{2n-3})+(c_1+c_{2n-2})t+(c_2+c_{2n-1})t^2 +(c_3+c_{2n})t^3}
\\ & & + \sum_{4\leq i\leq 2n-6} c_it^i   - c_{2n-5}
\Bigl(\sum_{i=0}^{\frac{2n-9}3} t^{3i+1} \Bigr)-
c_{2n-4} \Bigl( \sum_{i=0}^{\frac{2n-9}3} t^{3i+2}\Bigr) \nonumber  \\
& & + c_{2n-5}t\cdot{\textstyle \frac{1-t^{2n-3}}{1-t^3}} +
c_{2n-4}t^2\cdot{\textstyle \frac{1-t^{2n-3}}{1-t^3}} +
c_{2n-3}(t^{2n-3}-1)   \nonumber \\ & & + c_{2n-2}t(t^{2n-3}-1) +
c_{2n-1}t^2(t^{2n-3}-1)+ c_{2n}t^3(t^{2n-3}-1)  \nonumber
\end{eqnarray} where the equality comes from
\begin{eqnarray*} t^{2n-5} = -\sum_{i=0}^{\frac{2n-9}3} t^{3i+1}
+t\cdot{\textstyle \frac{1-t^{2n-3}}{1-t^3}}\ \ {\rm and}\ \
t^{2n-4} = -\sum_{i=0}^{\frac{2n-9}3} t^{3i+2} +
t^2\cdot{\textstyle \frac{1-t^{2n-3}}{1-t^3}} \end{eqnarray*}
since $\frac{1-t^{2n-3}}{1-t^3}=\displaystyle
\sum_{i=0}^{\frac{2n-6}3}t^{3i}$. Therefore, $p(X^{[n]};t)$ modulo
$\frac{1-t^{2n-3}}{1-t^3}$ is congruent to
\begin{eqnarray*} \lefteqn{R(t)=(c_0+c_{2n-3})+(c_1+c_{2n-2})t+(c_2+c_{2n-1})t^2
+(c_3+c_{2n})t^3} \\ & & + \displaystyle\sum_{4\leq i\leq 2n-6}
c_it^i - c_{2n-5} \Bigl(\displaystyle\sum_{i=0}^{\frac{2n-9}3}
t^{3i+1} \Bigr)- c_{2n-4} \Bigl(
\displaystyle\sum_{i=0}^{\frac{2n-9}3}
t^{3i+2}\Bigr).\end{eqnarray*} Now $R(t)$ is divisible by
$\frac{1-t^{2n-3}}{1-t^3}=\displaystyle
\sum_{i=0}^{\frac{2n-6}3}t^{3i}$ if and only if $R(t)$ is a scalar
multiple of $\displaystyle \sum_{i=0}^{\frac{2n-6}3}t^{3i}$
because $R(t)$ is of degree $\le 2n-6$. Thus the coefficient of
$R(t)$ with respect to $t^2$ should be 0 i.e.
$c_2+c_{2n-1}-c_{2n-4}=0$. However,
$c_2+c_{2n-1}-c_{2n-4}=c_2+c_1-c_4$ is not zero by Table
\ref{table:list of ci}. This proves Proposition \ref{prop:stringy
E-function test} for the case where 3 divides $n$ and $n\neq 3$.
So the proof of Proposition \ref{prop:stringy
E-function test} is completed for any $n\geq 3$. \qed\\

\begin{remark}
In case of smooth projective curves, we remark that the stringy
E-function of the moduli space of rank 2 bundles is explicitly
computed (\cite{kiem} and \cite{KL}). We were not able to compute
the stringy E-function of $M_{2n}$ precisely, because we do not
know how to compute the Hodge-Deligne polynomial $E(M^s_{2n};u,v)$
of the locus $M^s_{2n}$ of stable sheaves.
\end{remark}

%%%%%%%%%%%%%%%%%%%%%%%%%%%%%%%%%%%%%%%%%%%%%%%%%%%%%%%%%%%%%%%%%%%
%%%%%%%%%%%%%%%%%%%%%%%%%%%%%%%%%%%%%%%%%%%%%%%%%%%%%%%%%%%%%%%%%%%

\section{Analysis of Kirwan's desingularization}
\label{sec: Proof of Lemma}

This section is devoted to the proof of Proposition
\ref{prop:analysis on exc}. All can be extracted from \cite{og97}
but we spell out the details for reader's convenience.

To begin with, note that for each $Z\in X^{[n]}$, the tangent
space $T_{X^{[n]},Z}$ of the Hilbert scheme $X^{[n]}$ is
canonically isomorphic to $\Ext^1(I_Z,I_Z)$ where $I_Z$ is the
ideal sheaf of the 0-dimensional closed subscheme $Z$. By the
Yoneda pairing map and Serre duality, we have a skew-symmetric
pairing $\omega:\Ext^1(I_Z,I_Z)\otimes \Ext^1(I_Z,I_Z) \to
\Ext^2(I_Z,I_Z)\cong \cc$, which gives us a symplectic form
$\omega$ on the tangent bundle $T_{X^{[n]}}$ by \cite{Muk84}
Theorem 0.1.
%We write $E_Z=\Ext^1(I_Z,I_Z)$ for simplicity.
%And we
%use the notation $\omega$ for $\omega_Z$ if the meaning is clear
%in the context.

Note that the Killing form on $sl(2)$ gives an isomorphism
$sl(2)^\vee\cong sl(2)$. Let $W=sl(2)^\vee\cong sl(2)\cong \cc^3$.
The adjoint action of $PGL(2)$ on $W$ gives us an identification
$SO(W)\cong PGL(2)$ (\cite{og97} \S1.5). For a symplectic vector
space $(V,\omega)$, let $\Hom^\omega(W,V)$ be the space of
homomorphisms from $W$ to $V$ whose image is isotropic. Let
$\Hom^\omega(W,T_{X^{[n]}}) $ be the bundle over $X^{[n]}$ whose
fiber over $Z\in X^{[n]}$ is $\Hom^\omega(W,T_{X^{[n]},Z})$. %Since
%$X^{[n]}$ is smooth projective, $T_{X^{[n]}}$ is a Zariski locally
%trivial bundle. By elementary linear algebra, we can furthermore
%find local trivializations so that the symplectic form $\omega$ is
%given by a constant skew-symmetric matrix. Therefore, the bundle
Clearly $\Hom^\omega(W,T_{X^{[n]}})$ is Zariski locally trivial
over $X^{[n]}$. Let $\Hom_k^\omega(W,T_{X^{[n]}})$ be the
subbundle of $\Hom^\omega(W,T_{X^{[n]}})$ of rank $\leq k$
elements in $\Hom^\omega(W,T_{X^{[n]}})$. Also let
$\Gr^\omega(3,T_{X^{[n]}})$ be the relative Grassmannian of
isotropic 3-dimensional subspaces in $T_{X^{[n]}}$ and let $\cB$
denote the tautological rank 3 bundle on
$\Gr^\omega(3,T_{X^{[n]}})$. Obviously these bundles are all
Zariski locally trivial as well.

Let $\pp\Hom^\omega(W,T_{X^{[n]}})$ (resp.
$\pp\Hom_k^\omega(W,T_{X^{[n]}})$) be the projectivization of
$\Hom^\omega(W,T_{X^{[n]}})$ (resp.
$\Hom_k^\omega(W,T_{X^{[n]}})$). Likewise, let $\pp\Hom(W,\cB)$
and $\pp\Hom_k(W,\cB)$ denote the projectivizations of the bundles
$\Hom(W,\cB)$ and $\Hom_k(W,\cB)$. Note that there are obvious
forgetful maps
\begin{eqnarray*}f:\pp\Hom(W,\cB)\to\pp\Hom^\omega(W,T_{X^{[n]}})\ \mbox{\rm
and}\\
f_k:\pp\Hom_k(W,\cB)\to\pp\Hom_k^\omega(W,T_{X^{[n]}})\end{eqnarray*}
Since the pull-back of the defining ideal of
$\pp\Hom_1^\omega(W,T_{X^{[n]}})$ is the ideal of
$\pp\Hom_1(W,\cB)$ (both are actually given by the determinants of
$2\times 2$ minor matrices), $f$ gives rise to a map between
blow-ups
$$\overline{f}:Bl_{\pp\Hom_1(W,\cB)}\pp\Hom(W,\cB)\to
Bl_{\pp\Hom_1^\omega(W,T_{X^{[n]}})}\pp\Hom^\omega(W,T_{X^{[n]}}).$$
Let us denote $Bl_{\pp\Hom_1(W,\cB)}\pp\Hom(W,\cB)$ by $Bl^\cB$
and
$Bl_{\pp\Hom_1^\omega(W,T_{X^{[n]}})}\pp\Hom^\omega(W,T_{X^{[n]}})$
by $Bl^T$. We denote the proper transform of $\pp\Hom_2(W,\cB)$ in
$Bl^\cB$ by $Bl_2^\cB$ and the proper transform of
$\pp\Hom_2^\omega(W,T_{X^{[n]}})$ by $Bl_2^T$. Since $Bl_2^\cB$ is
a Cartier divisor which is mapped onto $Bl_2^T$ and the pull-back
of the defining ideal of $Bl_2^T$ is the ideal sheaf of
$Bl_2^\cB$, $\overline{f}$ lifts to
\begin{equation}\label{eq4.-2}
\hat{f}:Bl^\cB \to
Bl_{Bl_2^T}Bl^T.\end{equation} By \cite{og97} \S3.1 IV, $\hat f$
is an isomorphism on each fiber over $X^{[n]}$, so in particular
$\hat f$ is bijective. %Since both spaces are compact Hausdorff,
%$\hat f$ is a homeomorphism and thus by the Riemann extension
%theorem $\hat f^{-1}$ is holomorphic.
Therefore,  $\hat f$ is an isomorphism by Zariski's main theorem.

Note that $\pp\Hom(W,\cB)\git SO(W)$ (resp. $\pp\Hom_k(W,\cB)\git
SO(W)$) is isomorphic to the space of conics $\pp(S^2\cB)$ (resp.
rank $\leq k$ conics $\pp(S^2_k\cB)$)  where the $SO(W)$-quotient
map is given by $[\alpha]\mapsto[\alpha\circ\alpha^t]$ where
$\alpha^t$ denotes the transpose of $\alpha\in \Hom(W,\cB)$
(\cite{og97} \S3.1).  Let $\hat \pp(S^2\cB) =
Bl_{\pp(S^2_1\cB)}\pp(S^2\cB)$ denote the blow-up along the locus
of rank 1 conics. Then $Bl^\cB\git SO(W)$ is canonically
isomorphic to $\hat \pp(S^2\cB)$ by \cite{k2} Lemma 3.11. Since
$\cB$ is Zariski locally trivial, so is $\hat \pp(S^2\cB)$ over
$\Gr^\omega(3,T_{X^{[n]}})$.

Now consider Simpson's construction of the moduli space $M_{2n}$
(\cite{og97} \S1.1). Let $Q$ be the closure of the set of
semistable points $Q^{ss}$ in the Quot-scheme whose quotient by
the natural $PGL(N)$ action is $M_{2n}$ for some even integer $N$.
Then $Q^{ss}$ parameterizes semistable torsion-free sheaves $F$
together with surjective homomorphisms $h:\cO^{\oplus N}\to F(k)$
which induces an isomorphism $\cc^N\cong H^0(F(k))$ and
$H^1(F(k))=0$. Let $\Omega_Q$ denote the subset of $Q^{ss}$ which
parameterizes sheaves of the form $I_Z\oplus I_Z$ for some $Z\in
X^{[n]}$. This is precisely the locus of closed orbits with
maximal dimensional stabilizers, isomorphic to $PGL(2)$ and the
quotient of $\Omega_Q$ by $PGL(N)$ is $X^{[n]}$.

We can give a more precise description of $\Omega_Q$ as follows.
Let $\cL\to X^{[n]}\times X$ be the universal rank 1 sheaf such
that $\cL|_{Z\times X}$ is isomorphic to the ideal sheaf $I_Z$. By
\cite{HL97} Theorem 10.2.1, the tangent bundle $T_{X^{[n]}}$ is in
fact isomorphic to $\cE xt^1_{X^{[n]}}(\cL,\cL)$. Let
$p:X^{[n]}\times X\to X^{[n]}$ be the projection onto the first
component. For $k\gg 0$, $p_*\cL(k)$ is a vector bundle of rank
$N/2$. Let
\begin{equation}\label{eq4.-1}q:\pp \mathrm{Isom}(\cc^N, p_*\cL(k)\oplus
p_*\cL(k))\to X^{[n]}\end{equation} be the $PGL(N)$-bundle over
$X^{[n]}$ whose fiber over $Z$ is $\pp \mathrm{Isom}(\cc^N,
H^0(I_Z(k)\oplus I_Z(k)))$. Note that the standard action of
$GL(N)$ on $\cc^N$ and the obvious action of $GL(2)$ on
$p_*\cL(k)\oplus p_*\cL(k)$   induce a $PGL(N)\times
PGL(2)$-action  on $\pp\mathrm{Isom}(\cc^N, p_*\cL(k)\oplus
p_*\cL(k))\to X^{[n]}$ .
\begin{lemma}\label{4.1} (1) $\Omega_Q\cong \pp
\mathrm{Isom}(\cc^N, p_*\cL(k)\oplus p_*\cL(k))\git SO(W).$\\
(2) Via the above isomorphism, the normal cone of $\Omega_Q$ in
$Q^{ss}$ is $$q^*\mathrm{Hom}^{\omega}(W,T_{X^{[n]}})\git SO(W)\to
\pp \mathrm{Isom}(\cc^N, p_*\cL(k)\oplus p_*\cL(k))\git SO(W)$$
whose fiber  over a point lying over $Z\in X^{[n]}$ is $\mathrm{Hom}^{\omega}(W,T_{X^{[n]},Z})$. %Here
%$SO(W)$ acts on $q^*\mathrm{Hom}^{\omega}(W,T_{X^{[n]}})$ in the
%way that $SO(W)$ acts on $W$ and the base $PGL(N)$-bundle
%\eqref{eq4.-1}.
\end{lemma}
\begin{proof}
(1) Let $\hat p:\pp \mathrm{Isom}(\cc^N, p_*\cL(k)\oplus
p_*\cL(k))\times X\to \pp \mathrm{Isom}(\cc^N, p_*\cL(k)\oplus
p_*\cL(k))$ be the obvious projection so that we have $q\circ \hat
p=p\circ (q\times 1_X)$. Let $H$ be the dual of the tautological
line bundle over $\pp \mathrm{Isom}(\cc^N, p_*\cL(k)\oplus
p_*\cL(k))$. There is a canonical isomorphism $\cO^{\oplus N}\cong
q^*(p_*\cL(k)\oplus p_*\cL(k))\otimes H$.
 This induces a surjective homomorphism
\begin{eqnarray*}\cO^{\oplus N}\to \hat{p}^*q^*(p_*\cL(k)\oplus
p_*\cL(k))\otimes H =(q\times 1)^*(p^*p_*\cL(k)\oplus
p^*p_*\cL(k))\otimes H\\ \to (q\times 1)^*(\cL(k)\oplus
\cL(k))\otimes H\end{eqnarray*} over $\pp \mathrm{Isom}(\cc^N,
p_*\cL(k)\oplus p_*\cL(k))\times X$. By the universal property of
the Quot-scheme, we get a morphism $\pp \mathrm{Isom}(\cc^N,
p_*\cL(k)\oplus p_*\cL(k))\to Q^{ss}$ whose image is clearly
contained in $\Omega_Q$. This map is $PGL(2)$-invariant and hence
we get a  morphism
\begin{equation}\label{eq4.00}\phi_\Omega:\pp \mathrm{Isom}(\cc^N,
p_*\cL(k)\oplus p_*\cL(k))\git SO(W)\to \Omega_Q.\end{equation} It
is easy to check that $\phi_\Omega$ is bijective. Since $\Omega_Q$
is smooth (\cite{og97} (1.5.1)), $\phi_\Omega$ is an isomorphism
by Zariski's main theorem.

(2) Let  $\cO^{\oplus N} \to \cE(k)$ denote the universal quotient
sheaf on $Q^{ss}\times X$ restricted to $\Omega_Q$ and let $\cF$
be the kernel of the twisted homomorphism $\cO^{\oplus N}(-k) \to
\cE$ so that we have an exact sequence
$$0\to \cF \to \cO^{\oplus N}(-k) \to \cE\to 0$$ over
$\Omega_Q\times X$. The induced long exact sequence gives us
\begin{equation}\label{eqn: rel long exact sequence}
\HHom_{\Omega_Q}(\cO^{\oplus N}(-k),\cE)\to
\HHom_{\Omega_Q}(\cF,\cE)\to \EExt^1_{\Omega_Q}(\cE,\cE)\to
\EExt^1_{\Omega_Q}(\cO^{\oplus N}(-k),\cE) \end{equation} Let
$\pi:\Omega_Q\times X\to \Omega_Q$ be the obvious projection. Note
that $\EExt^1_{\Omega_Q}(\cO^{\oplus
N}(-k),\cE)=R^1\pi_*(\cE(k))^{\oplus N}=0$ and that
$\HHom_{\Omega_Q}(\cO^{\oplus N}(-k),\cE)\cong
\HHom_{\Omega_Q}(\cO^{\oplus N},\cE(k))$ is a vector bundle over
$\Omega_Q$ whose fiber is $gl(N)$ because $\cO_X^{\oplus N}\cong
H^0(E(k))$ for any $[\cO_X^{\oplus N}\to E(k)]\in Q^{ss}$. Let
$T^*_{Q^{ss}}, T^*_{\Omega_Q}$ be cotangent sheaves over $Q^{ss}$
and ${\Omega_Q}$ respectively.
%By the slice theorem (\cite{og97}
%(1.5.11)-(1.5.12)) we can check  that the following canonical
%sequence is exact (\cite{Hart} Theorem 8.17)
%$$0\to \cI_{\Omega_Q}/\cI_{\Omega_Q}^2\to
%T^*_{Q^{ss}}|_{\Omega_Q}\to T^*_{\Omega_Q}\to 0$$
By a famous result of Grothendieck
% Grothendieck's generalization of tangent spaces %%%%%%%%%%%%%%%%%%%%%%%%%%%%%%
(\cite{Gr95} \S5) we know
$$(T^*_{Q^{ss}}|_{\Omega_Q})^\vee\cong
\HHom_{\Omega_Q}(\cF,\cE)$$ which contains the tangent bundle of
$\Omega_Q$ as a subbundle. So the first homomorphism in
\eqref{eqn: rel long exact sequence} is the tangent map of the
group action of $PGL(N)$\footnote{In fact the term prior to the
first term of \eqref{eqn: rel long exact sequence} is
$\HHom_{\Omega_Q}(\cE,\cE)$ which contains $\cO$ obviously and the
quotient of $\HHom_{\Omega_Q}(\cO^{\oplus N}(-k),\cE)$ by $\cO$ is
a vector bundle whose fiber is the Lie algebra of $PGL(N)$. } on
$\Omega_Q$ and the second homomorphism is the Kodaira-Spencer map.

%It is easy to check that over each point $x\in \Omega_Q$ the image
%of $(T^*_{\Omega_Q})^\vee|_x\to
%(T^*_{Q^{ss}}|_{\Omega_Q})^\vee|_x$ is the image of
%$\Hom(\cO^{\oplus N}_X(-k),\cE|_x)\to \Hom(\cF|_x,\cE|_x)$ by
%\eqref{eqn: rel long exact sequence}. Therefore we obtain
%$(T^*_{Q^{ss}}|_{\Omega_Q})^\vee/(T^*_{\Omega_Q})^\vee\cong
%\EExt^1_{\Omega_Q}(\cE,\cE)$.

Via the isomorphism $\phi_\Omega$
\eqref{eq4.00}, we have a map
$$\delta:\pp
\mathrm{Isom}(\cc^N, p_*\cL(k)\oplus p_*\cL(k))\to \pp
\mathrm{Isom}(\cc^N, p_*\cL(k)\oplus p_*\cL(k))\git SO(W)\cong
\Omega_Q.$$ From the proof of (1) above, the pull-back of $\cE$ by
 $\delta\times 1$ is isomorphic to $(q\times 1)^*(\cL(k)\oplus
\cL(k))\otimes H$ and thus the vector bundle $\delta^*\cE
xt^1_{\Omega_Q}(\cE,\cE)$ is isomorphic to
$$q^*\cE xt^1_{X^{[n]}}(\cL, \cL)\otimes gl(2)\cong
q^*T_{X^{[n]}}\otimes gl(2).$$ The pull-back of the tangent sheaf
of $X^{[n]}$ sits in $q^*T_{X^{[n]}}\otimes gl(2)$ as
$q^*T_{X^{[n]}}\otimes \left(\begin{matrix}1&0\\
0&1\end{matrix}\right)$. % This is nothing but the pull-back by $\delta$
%of the image of $T_{\Omega_Q}$ by the Kodaira-Spencer map since when restricted to
%any closed point of $\Omega_Q$,
%the Kodaira-Spencer map maps $T_{\Omega_Q,x}$ onto $\Ext^1_{X^{[n]}}(\cL,\cL)\subset
%\Ext^1_{X^{[n]}}(\cL,\cL) \otimes gl(2)$.
Hence the pull-back by $\delta$ of the
normal bundle of $\Omega_Q$ (in the sense of \cite{og97} \S1.3) is
isomorphic to
$$q^*T_{X^{[n]}}\otimes sl(2)\cong
q^*\mathrm{Hom}(W,T_{X^{[n]}}).$$ By \cite{og97} (1.5.10), the
normal cone is  fiberwisely the same as the Hessian cone. (See
\cite{og97} \S1.3 for more details on the Hessian cone.) Since the
normal cone is contained in the Hessian cone, the normal cone is
equal to the Hessian cone which is the inverse image of zero by
the Yoneda square map $\Upsilon:\cE xt^1_{\Omega_Q}(\cE,\cE)\to
\cE xt^2_{\Omega_Q}(\cE,\cE)$. It is elementary to see that
$\delta^*\Upsilon^{-1}(0)$ is precisely
$q^*\mathrm{Hom}^\omega(W,T_{X^{[n]}}).$ Since $SO(W)$ acts freely
we obtain (2). See \cite{og97} (1.5.1) for a description of the
normal cone at each point.
\end{proof}

Let $\Sigma_Q$ denote the subset of $Q^{ss}$ whose sheaves are of
the form $I_Z\oplus I_W$ for some $Z,W\in X^{[n]}$. Then
$\Sigma_Q-\Omega_Q$ is precisely the locus of points in $Q^{ss}$
whose stabilizer is isomorphic to $\cc^*$. Let $\pi_R:R\to Q^{ss}$
be the blow-up of $Q^{ss}$ along $\Omega_Q$ and let $\Omega_R$
denote the exceptional divisor. By the above lemma,we have
\begin{equation}\label{eq4.0}\Omega_R\cong q^*\pp\mathrm{Hom}^{\omega}(W,T_{X^{[n]}})\git
SO(W).\end{equation} The following lemma is an easy exercise.
 \begin{lemma}\label{4.2} (1) The locus of points
in $\pp\mathrm{Hom}^\omega(W,T_{X^{[n]},Z})^{ss}$ whose stabilizer
is 1-dimensional by the action of $SO(W)$ is precisely $\pp
\mathrm{Hom}^\omega _1(W,T_{X^{[n]},Z})^{ss}$.\\ (2) The locus of
nontrivial stabilizers is $\pp \mathrm{Hom}^\omega
_2(W,T_{X^{[n]},Z})^{ss}$.\end{lemma} Let
\begin{equation}\label{eq4.1}\Delta_R=q^*\pp\mathrm{Hom}^{\omega}_2(W,T_{X^{[n]}})\git
SO(W).\end{equation} Let $\Sigma_R$ be the proper transform of
$\Sigma_Q$. Then $\Sigma_R^{ss}$ is precisely the locus of points
in $R^{ss}$ with 1-dimensional stabilizers by \cite{k2}. Therefore
we have the following from Lemma \ref{4.2}.
\begin{corollary}\label{4.3}
$\Sigma_R^{ss}\cap
\Omega_R=q^*\pp\mathrm{Hom}^{\omega}_1(W,T_{X^{[n]}})^{ss}\git
SO(W).$\end{corollary}

 We have an explicit description of
$\Sigma_R^{ss}$ from \cite{og97} \S1.7 III as follows. Let
$$\beta:\mathcal{X}^{[n]}\to X^{[n]}\times X^{[n]}$$
be the blow-up along the diagonal and let
$\mathcal{X}^{[n]}_0=X^{[n]}\times X^{[n]}-\mathbf{\Delta}$ where
$\mathbf{\Delta}$ is the diagonal. Let $\cL_1$ (resp. $\cL_2$) be
the pull-back to $\mathcal{X}^{[n]}\times X$ of the universal
sheaf $\cL\to X^{[n]}\times X$ by  $p_{13}\circ (\beta\times 1)$
(resp. $p_{23}\circ (\beta\times 1)$) where $p_{ij}$ is the
projection onto the first (resp. second) and third components. Let
$p:\mathcal{X}^{[n]}\times X\to \mathcal{X}^{[n]}$ be the
projection onto the first component. Then for   $k\gg0$,
$p_*\cL_1(k)\oplus p_*\cL_2(k)$ is a vector bundle of rank $N$.
Let
$$q:\pp\mathrm{Isom}(\cc^N,p_*\cL_1(k)\oplus p_*\cL_2(k))\to \mathcal{X}^{[n]}$$
be the $PGL(N)$-bundle. There is an action of $O(2)$ on
$\pp\mathrm{Isom}(\cc^N,p_*\cL_1(k)\oplus p_*\cL_2(k))$. We quote
 \cite{og97} (1.7.10) and (1.7.1).
\begin{lemma}\label{4.4}
(1) $\Sigma_R^{ss}\cong\pp\mathrm{Isom}(\cc^N,p_*\cL_1(k)\oplus
p_*\cL_2(k))\git O(2)$\\
(2) The normal cone of $\Sigma_R^{ss}$ in $R^{ss}$ is a locally
trivial bundle over $\Sigma_R^{ss}$ with fiber the cone over a
smooth quadric in $\pp^{4n-5}$.
\end{lemma}
In fact we can give a more explicit description of the normal cone
when restricted to $\Sigma_R^0:=\Sigma_R^{ss}-\Omega_R$. Similarly
as in the proof of Lemma \ref{4.1}, the normal vector bundle to
$\Sigma_R^0$ is isomorphic to the vector bundle (of rank $4n-4$)
\begin{equation}\label{eq4.2} q^*[\cE
xt^1_{\mathcal{X}^{[n]}_0}(\cL_1,\cL_2)\oplus \cE
xt^1_{\mathcal{X}^{[n]}_0}(\cL_2,\cL_1)]\git O(2)\end{equation}
over $\pp\mathrm{Isom}(\cc^N,p_*\cL_1(k)\oplus p_*\cL_2(k))\git
O(2)$ where $O(2)$ acts as follows: if we realize $O(2)$ as the
subgroup of $PGL(2)$
generated by $$SO(2)=\{\theta_\alpha=\left(\begin{matrix}\alpha&0\\
0&\alpha^{-1}\end{matrix}\right)\}/\{\pm Id\},\qquad
\tau=\left(\begin{matrix} 0&1\\1&0\end{matrix}\right)$$
$\theta_\alpha$ multiplies $\alpha$ (resp. $\alpha^{-1}$) to
$\cL_1$ (resp. $\cL_2$) and $\tau$ interchanges $\cL_1$ and
$\cL_2$ by the induced action on $\mathcal{X}^{[n]}$ of
interchanging the first and second factors of $X^{[n]}\times
X^{[n]}$. The normal cone is the inverse image
$q^*\Upsilon^{-1}(0)$ of zero in terms of the Yoneda pairing
\begin{equation}\label{eq4.3}\Upsilon:\cE
xt^1_{\mathcal{X}^{[n]}_0}(\cL_1,\cL_2)\oplus \cE
xt^1_{\mathcal{X}^{[n]}_0}(\cL_2,\cL_1)\to \cE
xt^2_{\mathcal{X}^{[n]}_0}(\cL_1,\cL_1).\end{equation}

Let $\pi_S:S\to R^{ss}$ denote the blow-up of $R^{ss}$ along
$\Sigma_R^{ss}$ and let $\Sigma_S$ be the exceptional divisor of
$\pi_S$ while $\Omega_S$ (resp. $\Delta_S$) denotes the proper
transform of $\Omega_R$ (resp. $\Delta_R$). By \eqref{eq4.3}, we
have \begin{equation}\label{eq4.4} {\Sigma_S}
|_{\pi_S^{-1}(\Sigma_R^0)}\cong q^*\pp \Upsilon^{-1}(0)\git
O(2)\subset q^*\pp[\cE
xt^1_{\mathcal{X}^{[n]}_0}(\cL_1,\cL_2)\oplus \cE
xt^1_{\mathcal{X}^{[n]}_0}(\cL_2,\cL_1)]\git O(2).
\end{equation}

 By
\cite{og97} (1.8.10), $S^s=S^{ss}$ and $S^s$ is smooth. The
quotient $S\git PGL(N)$ has only $\zz_2$-quotient singularities
along $\Delta_S\git PGL(N)$. Let $\pi_T:T\to S^{s}$ be the blow-up
of $S^{s}$ along $\Delta_S^{s}$. Then $T\git PGL(N)$ is
nonsingular and this is Kirwan's desingularization
$\rho:\hM_{2n}\to M_{2n}$.

Let $\Omega_T$ and $\Sigma_T$ denote the proper transforms of
$\Omega_S$ and $\Sigma_S$ respectively. Let $\Delta_T$ be the
exceptional divisor of $\pi_T$. Their quotients $\Omega_T\git
PGL(N)$, $\Sigma_T\git PGL(N)$ and $\Delta_T\git PGL(N)$ are
denoted by $D_1=\hat\Omega$, $D_2=\hat\Sigma$ and $D_3=\hat\Delta$
respectively.

With this preparation, we now embark on the proof of Proposition
\ref{prop:analysis on exc}.

\vspace{.5cm}\noindent \textbf{Proof of (1).}  This is just
\cite{og97} (3.0.1). More precisely, by \eqref{eq4.0} and
Corollary \ref{4.3}, $\Omega_S$ is the blow-up of
$$q^*\pp\mathrm{Hom}^{\omega}(W,T_{X^{[n]}})\git SO(W)\text{ along
}q^*\pp\mathrm{Hom}^{\omega}_1(W,T_{X^{[n]}})\git SO(W).$$ By
\eqref{eq4.1}, $\Omega_T$ is the blow-up of $\Omega_S$ along the
proper transform of
$$q^*\pp\mathrm{Hom}^{\omega}_2(W,T_{X^{[n]}})\git SO(W)$$ and
$D_1=\hat\Omega$ is the quotient of $\Omega_T$ by the action of
$PGL(N)$. Since the action of $PGL(N)$ commutes with the action of
$SO(W)$, $D_1$ is in fact the quotient by $SO(W)\times PGL(N)$ of
the variety obtained from
$q^*\pp\mathrm{Hom}^{\omega}(W,T_{X^{[n]}})$ by two blow-ups. So
$D_1$ is also the consequence of taking the quotient by $PGL(N)$
first and then the quotient by $SO(W)$ second. Since $q$
\eqref{eq4.-1} is a principal $PGL(N)$ bundle, the result of the
first quotient is just $Bl_{Bl_2^T}Bl^T$ in \eqref{eq4.-2} which
is isomorphic to $Bl^\cB$. If we take further the quotient by
$SO(W)$, then as discussed above the result is $D_1=\hat\pp
(S^2\cB)$.

\vspace{.5cm}\noindent \textbf{Proof of (2).} We use Lemma
\ref{4.4}, \eqref{eq4.2}, and \eqref{eq4.4}. Note that
$\Sigma_R^0$ does not intersect with $\Omega_R$ and $\Delta_R$.
Hence $D_2^0$ is the quotient of $q^*\pp \Upsilon^{-1}(0)\git
O(2)$ which is a subset of $q^*\pp[\cE
xt^1_{\mathcal{X}^{[n]}_0}(\cL_1,\cL_2)\oplus \cE
xt^1_{\mathcal{X}^{[n]}_0}(\cL_2,\cL_1)]\git O(2)$, by the action
of $PGL(N)$. The above are bundles over the restriction of
$$\pp\mathrm{Isom}(\cc^N,p_*\cL_1(k)\oplus p_*\cL_2(k))\git O(2)$$
to the complement $\mathcal{X}^{[n]}_0$ of the diagonal
$\mathbf{\Delta}$ in $X^{[n]}\times X^{[n]}$. As in the proof of
(1), observe that $D_2^0$ is in fact the quotient of
$q^*\pp\Upsilon^{-1}(0)$ by the action of $PGL(N)\times O(2)$
since the actions commute. So we can first take the quotient by
the action of $PGL(N)$, then by the action of $SO(2)$, and finally
by the action of $\zz_2=O(2)/SO(2)$. Since
$\pp\mathrm{Isom}(\cc^N,p_*\cL_1(k)\oplus p_*\cL_2(k))$ is a
principal $PGL(N)$-bundle, the quotient by $PGL(N)$ gives us
$$\pp\Upsilon^{-1}(0)\subset \pp[\cE
xt^1_{\mathcal{X}^{[n]}_0}(\cL_1,\cL_2)\oplus \cE
xt^1_{\mathcal{X}^{[n]}_0}(\cL_2,\cL_1)]$$ over
$\mathcal{X}^{[n]}_0$. The algebraic vector bundles $\cE
xt^1_{\mathcal{X}^{[n]}_0}(\cL_1,\cL_2)$ and $\cE
xt^1_{\mathcal{X}^{[n]}_0}(\cL_2,\cL_1)$ are certainly Zariski
locally trivial and in fact these bundles are dual to each other
by the Yoneda pairing $\Upsilon$ which is non-degenerate (possibly
after tensoring with a line bundle). In particular,
$\Upsilon^{-1}(0)$ is Zariski locally trivial.

Next we take the quotient by the action of $SO(2)\cong \cc^*$.
This action is trivial on the base $\mathcal{X}^{[n]}_0$ and
$SO(2)$ acts on the fibers. Hence $\pp\Upsilon^{-1}(0)/SO(2)$ is a
Zariski locally trivial subbundle of
$$\pp[\cE xt^1_{\mathcal{X}^{[n]}_0}(\cL_1,\cL_2)\oplus \cE
xt^1_{\mathcal{X}^{[n]}_0}(\cL_2,\cL_1)]\git \cc^*\cong \pp\cE
xt^1_{\mathcal{X}^{[n]}_0}(\cL_1,\cL_2)\times_{\mathcal{X}^{[n]}_0}
\pp\cE xt^1_{\mathcal{X}^{[n]}_0}(\cL_2,\cL_1)$$ over
$\mathcal{X}^{[n]}_0$ given by the incidence relations in terms of
the identification $$\pp\cE
xt^1_{\mathcal{X}^{[n]}_0}(\cL_1,\cL_2)\cong \pp\cE
xt^1_{\mathcal{X}^{[n]}_0}(\cL_2,\cL_1)^\vee.$$ Finally, $D_2^0$
is the $\zz_2$-quotient of $\pp\Upsilon^{-1}(0)/SO(2)$.

\vspace{.5cm}\noindent \textbf{Proof of (3).} By \cite{og97}
(1.7.10), the intersection of $\Sigma_R^{ss}$ and $\Omega_R$ is
smooth. By Corollary \ref{4.3}, $\Delta_S$ is the blow-up of
$q^*\pp\mathrm{Hom}^{\omega}_2(W,T_{X^{[n]}})\git SO(W)$ along
$q^*\pp\mathrm{Hom}^{\omega}_1(W,T_{X^{[n]}})\git SO(W)$. Hence
$\Delta_S\git PGL(N)$ is the quotient of  $$
Bl_{q^*\pp\mathrm{Hom}^{\omega}_1
(W,T_{X^{[n]}})}q^*\pp\mathrm{Hom}^{\omega}_2(W,T_{X^{[n]}})$$ by
the action of  $SO(W)\times PGL(N)$. By taking the quotient by the
action of $PGL(N)$ we get
$$ Bl_{ \pp\mathrm{Hom}^{\omega}_1
(W,T_{X^{[n]}})} \pp\mathrm{Hom}^{\omega}_2(W,T_{X^{[n]}})$$since
$q$ is a principal $PGL(N)$-bundle. Next we take the quotient  by
the action of $SO(W)$. Let $\mathrm{Gr}^\omega(2,T_{X^{[n]}})$ be
the relative Grassmannian of isotropic 2-dimensional subspaces in
$T_{X^{[n]}}$ and let $\mathcal A$ be the tautological rank 2
bundle on $\mathrm{Gr}^\omega(2,T_{X^{[n]}})$. We claim
\begin{equation}\label{}Bl_{\pp \Hom_1^\omega(W,T_{X^{[n]}})}
{\pp \Hom_2^\omega(W,T_{X^{[n]}})}\git SO(W)\simeq \pp(S^2\cA)
\end{equation} which is a $\pp^{2}$-bundle
over a  $\Gr^\omega (2,2n)$-bundle over $X^{[n]}$. It is obvious
that the bundles are Zariski locally trivial.

There are forgetful maps \begin{eqnarray*}f:\pp\Hom(W,\cA)\to
\pp\Hom_2^\omega(W,T_{X^{[n]}}) \\ f_1:\pp\Hom_1(W,\cA)\to
\pp\Hom_1^\omega(W,T_{X^{[n]}}) \end{eqnarray*} where the
subscript 1 denotes the locus of rank $\leq1$ homomorphisms.
Because the ideal of $\pp\Hom^\omega_1(W,T_{[n]})$ pulls back to
the ideal of $\pp\Hom_1(W,\cA)$, $f$ lifts to $$ \hat
f:Bl_{\pp\Hom_1(W,\cA)} \pp\Hom(W,\cA)\to
Bl_{\pp\Hom_1^\omega(W,T_{X^{[n]}})}
{\pp\Hom_2^\omega(W,T_{X^{[n]}})}$$ This map is bijective
(\cite{og97} (3.5.1)) and hence $\hat f$ is an isomorphism by
Zariski's main theorem because the varieties are smooth. Now
observe that the quotient $\pp\Hom(W,\cA)\git SO(W)$ is
$\pp(S^2\cA)$ where the quotient map is given by $[\alpha]\mapsto
[\alpha\circ\alpha^t]$. Hence $\Delta_S\git PGL(N) $ is the
blow-up of $\pp\Hom(W,\cA)\git SO(W)\cong\pp(S^2\cA)$ along the
locus of rank 1 quadratic forms $\pp(S^2_1\cA)$ (\cite{k2} Lemma
3.11) which is a Cartier divisor. So we proved that
$$ \Delta_S\git PGL(N)\cong \pp(S^2\cA).$$

Finally $S\git PGL(N)$ is singular only along
$\Delta_S\git PGL(N)$ and the singularities are $\cc^{2n-3}/\{\pm
1\}$ by Luna's slice theorem \cite{og97} (1.2.1). Since $D_3$ is
the exceptional divisor of the blow-up of $S\git PGL(N)$ along
$\Delta_S\git PGL(N)$, we conclude that $D_3$ is a
$\pp^{2n-4}$-bundle over $\pp (S^2\mathcal A)$.

\vspace{.5cm}\noindent \textbf{Proof of (4).} By Corollary
\ref{4.3}, $\Sigma_S^s\cap\Omega_S$ is the exceptional divisor of
the blow-up
$Bl_{q^*\pp\mathrm{Hom}^{\omega}_1(W,T_{X^{[n]}})}q^*\pp\mathrm{Hom}^{\omega}(W,T_{X^{[n]}})\git
SO(W) $ and $\Sigma_T^s\cap \Omega_T$ is now the blow-up of the
exceptional divisor along the proper transform of
$q^*\pp\mathrm{Hom}^{\omega}_2(W,T_{X^{[n]}})\git SO(W)$. Using
the isomorphism \eqref{eq4.-2}, this is the exceptional divisor of
$$q^*Bl_{\pp (S^2_1\cB)}\pp (S^2\cB)\to q^*\pp (S^2\cB)$$ over
$\mathrm{Gr}^\omega(3,T_{X^{[n]}})$. Since $q$ is a principal
$PGL(N)$-bundle, $D_1\cap D_2=\Sigma_T^s\cap \Omega_T\git PGL(N)$
is the exceptional divisor of the blow-up $Bl_{\pp (S^2_1\cB)}\pp
(S^2\cB)$.   Because the exceptional divisor is a Zariski locally
trivial $\pp^2 $-bundle over  $\pp (S^2_1\cB)$ and $\pp
(S^2_1\cB)$ itself is a Zariski locally trivial $\pp^2$-bundle
over $\mathrm{Gr}^\omega(3,T_{X^{[n]}})$, we proved (4).

\vspace{.5cm}\noindent \textbf{Proof of (5).} From the above proof
of (3) it follows immediately that $\Sigma_S^s\cap \Delta_S\git
PGL(N)$ is $\pp (S^2_1\cA)$ and $D_2\cap D_3$ is a $\pp^{2n-4}$
bundle over $\pp (S^2_1\cA)$ which is Zariski locally trivial.

\vspace{.5cm}\noindent \textbf{Proof of (6).} As in the above
proof of (4), we start with \eqref{eq4.1} and use the isomorphism
\eqref{eq4.-2} to see that $D_1\cap D_3$ is the proper transform
of $\pp (S^2_2\cB)$ in the blow-up $Bl_{\pp (S^2_1\cB)}\pp
(S^2\cB)$. This is a Zariski locally trivial $\pp^2$-bundle over a
Zariski locally trivially $\pp^2$-bundle over
$\mathrm{Gr}^\omega(3,T_{X^{[n]}})$.

\vspace{.5cm}\noindent \textbf{Proof of (7).} This follows
immediately from the proof of (4) and (6).

\vspace{.5cm}
 From the above descriptions, it is clear that $D_i$
($i=1,2,3$) are normal crossing smooth divisors. \qed

%%%%%%%%%%%%%%%%%%%%%%%%%%%%%%%%%%%%%%%%%%%%%%%%%%%%%%%%%%%%%%%%%%
%%%%%%%%%%%%%%%%%%%%%%%%%%%%%%%%%%%%%%%%%%%%%%%%%%%%%%%%%%%%%%%%%%

\section{Hodge-Deligne polynomial of $D_2^0 $}
\label{sec: Computation of E-poly of D_0^2}

In this section we prove Lemma \ref{lem: Hodge Deligne poly of
D02}. Recall $$I_{2n-3}=\{((x_i),(y_j))\in \pp^{2n-3}\times
\pp^{2n-3}\,|\, \sum_{i=0}^{2n-3} x_iy_i=0\}.$$ It is elementary
(\cite{GH78} p. 606) to see that
$$H^*(I_{2n-3};\qq)\cong \qq[a,b]/\langle a^{2n-2}, b^{2n-2},
a^{2n-3}+a^{2n-4}b+a^{2n-5}b^2+\cdots+b^{2n-3} \rangle$$ where $a$
(resp. $b$) is the pull-back of the first Chern class of the
tautological line bundle of the first (resp. second) $\pp^{2n-3}$.
The $\zz_2$-action interchanges $a$ and $b$ and the invariant
subspace of $H^*(I_{2n-3};\qq)$ is generated by classes of the
form $a^ib^j+a^jb^i$. As a vector space $H^*(I_{2n-3};\qq)$ is
\begin{equation}\label{eq5.0}\qq\text{-span}\{a^ib^j\,|\, 0\le i\le 2n-3, 0\le j\le
2n-4\}\end{equation} while the invariant subspace is
$$\qq\text{-span}\{a^ib^j+a^jb^i\,|\, 0\le i\le j\le 2n-4\}.$$
The index set $\{(i,j)\,|\, 0\le i\le j\le 2n-4\}$ is mapped to
its complement in $\{(i,j)\,|\, 0\le i\le 2n-3, 0\le j\le 2n-4\}$
by the map $(i,j)\mapsto (j+1,i)$. This immediately implies that
\begin{equation}\label{eq5.3}
P(I_{2n-3};z)=(1+z^2)P^+(I_{2n-3};z)\end{equation}
 By \eqref{eq5.0} or the observation that $I_{2n-3}$ is the
Zariski locally trivial $\pp^{2n-4}$-bundle over $\pp^{2n-3}$, we
have
\begin{equation}\label{eq5.2}P(I_{2n-3};z)=\frac{1-(z^2)^{2n-2}}{1-z^2}\cdot
\frac{1-(z^2)^{2n-3}}{1-z^2}.\end{equation} Because $1+z^2$
divides $\frac{1-(z^2)^{2n-2}}{1-z^2}$,
$\frac{1-(z^2)^{2n-3}}{1-z^2}$ also divides $P^+(I_{2n-3};z)$.
Therefore, (\ref{eqn: E D02 is divisible by some Q}) is a direct
consequence of \eqref{eqn: compute D02} since $P(X^{[n]};z)$ has
no odd degree terms by \eqref{eqn:Betti for X[n]}.

Now let us prove \eqref{eqn: compute D02}. Let
$$\psi:\widetilde{D}_2^0:=\pp\Upsilon^{-1}(0)/SO(2)\to
\mathcal{X}^{[n]}_0=X^{[n]}\times X^{[n]}-\mathbf{\Delta}$$ be the
Zariski locally trivial $I_{2n-3}$-bundle in the proof of
Proposition \ref{prop:analysis on exc} (2) in \S4. Recall that
$D_2^0=\widetilde{D}_2^0/\zz_2$. We have seen in the proof of
Proposition \ref{prop:analysis on exc} (2) in \S4 that there is a
$\zz_2$-equivariant embedding
$$\imath:\widetilde{D}_2^0\hookrightarrow
\pp\cE
xt^1_{\mathcal{X}^{[n]}_0}(\cL_1,\cL_2)\times_{\mathcal{X}^{[n]}_0}
\pp\cE xt^1_{\mathcal{X}^{[n]}_0}(\cL_2,\cL_1)$$ where the
$\zz_2$-action interchanges $\cL_1$ and $\cL_2$.

Let $\lambda$ (resp. $\eta$) be the pull-back to
$\widetilde{D}_2^0$ of the first Chern class of the tautological
line bundle over $\pp\cE xt^1_{\mathcal{X}^{[n]}_0}(\cL_1,\cL_2)$
(resp. $\pp\cE xt^1_{\mathcal{X}^{[n]}_0}(\cL_2,\cL_1)$). By
definition, $\lambda$ and $\eta$ restrict to $a$ and $b$
respectively. The $\zz_2$-action interchanges $\lambda$ and
$\eta$. By the Leray-Hirsch theorem\footnote{The Leray-Hirsch
theorem in \cite{V02I} p.182 is stated for ordinary cohomology but
the statement holds also for compact support cohomology. See the
proof in \cite{V02I} p.195} we have an isomorphism
\begin{equation}\label{eq5.5} H^*_c(\widetilde{D}_2^0)\ \
%\cong \bigoplus_{0\le i\le 2n-3, 0\le j\le 2n-4 }
%H^*_c(\mathcal{X}^{[n]}_0)\lambda^i\eta^j\
\cong \ H^*_c(\mathcal{X}^{[n]}_0)\otimes
H^*(I_{2n-3}).\end{equation} As the pull-back and the cup product
preserve mixed Hodge structure, \eqref{eq5.5} determines the mixed
Hodge structure of $H^*_c(\widetilde{D}_2^0)$. The
$\zz_2$-invariant part is
\begin{equation}\label{eq5.6}
H^*_c(\widetilde{D}_2^0)^+\cong \left(
H^*_c(\mathcal{X}^{[n]}_0)^+\otimes H^*(I_{2n-3})^+\right) \oplus
\left( H^*_c(\mathcal{X}^{[n]}_0)^-\otimes H^*(I_{2n-3})^-\right)
\end{equation} where the superscript $\pm$ denotes the $\pm
1$-eigenspace of the $\zz_2$-action. Because $H^*_c(D_2^0)\cong
H^*_c(\widetilde{D}_2^0/\zz_2)\cong H^*_c(\widetilde{D}_2^0)^+$
(\cite{Gr57} Theorem 5.3.1 and Proposition 5.2.3), $E(D_2^0;u,v)$
is equal to
\begin{equation}\label{eq5.7}
 E^+(\widetilde{D}_2^0;u,v)=E^+(\mathcal{X}^{[n]}_0;u,v)E^+(I_{2n-3};u,v)
 +E^-(\mathcal{X}^{[n]}_0;u,v)E^-(I_{2n-3};u,v).
 \end{equation}
where $E^\pm(Y;u,v)=\sum_{p,q}\sum_{k\geq0} (-1)^k
h^{p,q}(H^k_c(Y)^\pm) u^pv^q$.

It is easy to see \begin{eqnarray*}P^+(X^{[n]}\times
X^{[n]};z)=\frac{P(X^{[n]};z)^2+P(X^{[n]};z^2)}2,\\
{P^-(X^{[n]}\times
X^{[n]};z)=\frac{P(X^{[n]};z)^2-P(X^{[n]};z^2)}2}\end{eqnarray*}
(Macdonald's formula). Since $X^{[n]}\times X^{[n]}$ is smooth
projective, we have
\begin{eqnarray*} E^+(X^{[n]}\times
X^{[n]};z,z)=\frac{P(X^{[n]};z)^2+P(X^{[n]};z^2)}2\\
E^-(X^{[n]}\times
X^{[n]};z,z)=\frac{P(X^{[n]};z)^2-P(X^{[n]};z^2)}2
\end{eqnarray*}
Now as $\mathcal{X}^{[n]}_0=X^{[n]}\times X^{[n]}-\mathbf{\Delta}$
and $\mathbf{\Delta}\cong X^{[n]}$ is $\zz_2$-invariant, by the
additive property of the E-polynomial we have
\begin{eqnarray*}
E^+(\mathcal{X}^{[n]}_0;z,z)=E^+(X^{[n]}\times X^{[n]};z,z)
-E(X^{[n]};z,z)\\ =
\frac{P(X^{[n]};z)^2+P(X^{[n]};z^2)}2-P(X^{[n]};z),\end{eqnarray*}
\begin{eqnarray*}E^-(\mathcal{X}^{[n]}_0;z,z)=E^-(X^{[n]}\times X^{[n]};z,z)\\
=\frac{P(X^{[n]};z)^2-P(X^{[n]};z^2)}2.\end{eqnarray*} The
equation \eqref{eqn: compute D02} is an immediate consequence of
the above equations and \eqref{eq5.7}. \qed

\bibliographystyle{amsplain}

\end{document}